\documentclass[reqno, 12pt]{amsart}
\usepackage{graphicx}

\newtheorem{thm}{Theorem}[section]

\theoremstyle{definition}

\theoremstyle{remark}

\numberwithin{equation}{section}
\begin{document}
\hfill{To appear in: \emph{Quality Technology and Quantitative
Management\\}

\title[Asymptotic Analysis of Loss Probabilities]{Asymptotic Analysis of Loss Probabilities in
$GI/M/m/n$ Queueing Systems as $n$ Increases to Infinity}%
\author{Vyacheslav M. Abramov}%
\address{School of Mathematical Sciences, Monash University,
Building 28M, Clayton campus, Clayton, Victoria 3800, Australia}%
\email{Vyacheslav.Abramov@sci.monash.edu.au}%

%\thanks{}%

%\date{}%
%\dedicatory{}%
%\commby{}%
% ----------------------------------------------------------------
\begin{abstract}
The paper studies asymptotic behavior of the loss probability for
the $GI/M/m/n$ queueing system as $n$ increases to infinity. The
approach of the paper is based on applications of classic results
of Tak\'acs (1967) and the Tauberian theorem with remainder of
Postnikov (1979-1980) associated with the recurrence relation of
convolution type. The main result of the paper is associated with
asymptotic behavior of the loss probability. Specifically it is
shown that in some cases (precisely described in the paper) where
the load of the system approaches 1 from the left and $n$
increases to infinity, the loss probability of the $GI/M/m/n$
queue becomes asymptotically independent of the parameter $m$.
\end{abstract}
\maketitle

\newpage

% ----------------------------------------------------------------
\section{\bf Introduction}
\noindent It is well-known that queueing systems with many servers
are well models for communication systems. Study of queues with
many servers and especially analysis of the loss probability have
a long history going back to the works of Erlang (see \cite{the
life of Erlang}), who in 1917 first gave fundamental results for
Markovian queueing systems, and to the works of Palm, Pollaczek
and other researchers (e.g. \cite{Palm (1943)}, \cite{Pollaczek
(1953)}, \cite{Cohen (1957)}, \cite{Sevastyanov (1957)},
\cite{Takacs (1957)}), who then developed the Erlang's results to
non-Markovian systems. Since then these results have been
developed in a large number of investigations, motivated by
growing development of modern telecommunication systems. Nowadays
the theory of loss queueing theory is very rich. There is a number
of different directions of the theory including management and
control, redundancy, analysis of retrials, impatient customers and
so on.

In the present paper we study the $GI/M/m/n$ queueing system in
which parameter $n$, the number of possible waiting places, is a
large value. This assumption is typical for real telecommunication
systems. Under this assumption the paper studies asymptotic
behavior of the loss probability, where the most significant
result seems to be related to the case where the load parameter
approaches 1 from the left. The case of a heavy load parameter is
the most interesting in practice. In the pre-project study and
system design stage an engineer is especially interested to know
about the behavior of the loss probability in heavy loaded
systems.

We consider the $GI/M/m/n$ queue, where $m>1$ is the number of
servers, and $n\geq 0$ is the admissible queue-length. Let $A(x)$
denote the probability distribution function of an interarrival
time, and let $\lambda$ be the reciprocal of the expected
interarrival time. For real $s\geq 0$ we denote
$\alpha(s)=\int_0^\infty\mbox{e}^{-sx}\mbox{d}A(x)$. The parameter
of the service time distribution is denoted $\mu$, and the load of
the system is $\varrho=\lambda/(m\mu)$. (In Theorem \ref{t2} the
parameter $\varrho$ is assumed to depend on $n$. However we do not
write this dependence explicitly, assuming that this dependence is
clear from the formulation of the aforementioned theorem.)

In the case of $GI/M/1/n$ queueing system for the stationary loss
probability $p_n$ we have the representation (Abramov
\cite{Abramov (2002)})
\begin{equation}\label{r1.1}
 p_n=\frac{1}{\pi_n},
\end{equation}
where the generating function $\Pi(z)$ of $\pi_j$, $j=0,1,\ldots$
is the following:
\begin{equation}\label{r1.2}
\Pi(z)=\sum_{j=0}^\infty\pi_jz^j=\frac{\alpha(\mu-\mu
z)}{\alpha(\mu-\mu z)-z}, \ \ |z|<\sigma,
\end{equation}
$\sigma$ is the least in absolute value root of the functional
equation $z=\alpha(\mu-\mu z)$ (the variable $z$ is assumed to be
real). It is well-known (e.g. Tak\'acs \cite{Takacs (1962)}), that
$\sigma$ belongs to the open interval (0,1) if $\varrho<1$, and it
is equal to 1 otherwise. Note, that another representation for the
loss probability $p_n$ is given in Miyazawa \cite{Miyazawa
(1990)}.

The value $\pi_n$ has the following meaning. This is the expected
number of arrivals up to the first loss of a customer arriving to
the stationary system. $\pi_n$ satisfies the recurrence relation
of convolution type
\begin{equation}\label{r1.3}
\pi_n=\sum_{i=0}^n r_i\pi_{n-i+1}, \ \ n=0,1,\ldots,
\end{equation}
where the initial value is $\pi_0=1$. Specifically,
\begin{equation}\label{r1.4}
r_i=\int_0^\infty\mbox{e}^{-\mu x}\frac{(\mu
x)^i}{i!}\mbox{d}A(x).
\end{equation}
As $n\to\infty$, the asymptotic behavior of the general recurrence
relation of convolution type,
\begin{equation}\label{r1.5}
Q_n=\sum_{i=0}^n f_iQ_{n-i+1},
\end{equation}
where $f_0>0$, $f_i\geq 0$ ($i\geq 1$) and $f_0+f_1+...=1$, has
been originally studied by Tak\'acs \cite{Takacs (1967)}, p. 22
and then developed by Postnikov \cite{Postnikov (1979)}, Section
25. For the readers' convenience, some of these results, necessary
for the purpose of the paper, are collected in the Appendix.

By exploiting (\ref{r1.3}), Abramov \cite{Abramov (2002)} studied
the asymptotic behavior of the loss probability $p_n$, as
$n\to\infty$. The analysis of Abramov \cite{Abramov (2002)} is
based on the mentioned results of Tak\'acs \cite{Takacs (1967)}
and Postnikov \cite{Postnikov (1979)}. For other applications of
the mentioned results of Tak\'acs \cite{Takacs (1967)} and
Postnikov \cite{Postnikov (1979)} see also Abramov \cite{Abramov
(1997)}, \cite{Abramov (2003)}, \cite{Abramov (2005)}.

The asymptotic analysis of $GI/M/m/n$ queueing system, as
$n\to\infty$, is a much more difficult problem than the same
problem for the $GI/M/1/n$ queue. Recently, Choi et al \cite{Choi
Kim Kim and Wee (2003)} and Kim and Choi \cite{Kim Choi (2003)}
obtained some new results related to the $GI/M/m/n$ and
$GI^X/M/m/n$ queues. In Choi et al \cite{Choi Kim Kim and Wee
(2003)} the exact estimation for the convergence rate of the
stationary $GI/M/m/n$ queue-length distribution to the stationary
queue-length distribution of the $GI/M/m$ queueing system, as
$n\to\infty$, is obtained.  In Kim and Choi \cite{Kim Choi (2003)}
detailed analysis of the loss probability of the $GI^X/M/m/n$ is
provided. The analysis of these two aforementioned papers is based
on the development of the earlier results of Choi and Kim
\cite{Choi and Kim (2000)} in a nontrivial fashion.

The analysis of Choi and Kim \cite{Choi and Kim (2000)} and Choi
et al \cite{Choi Kim and Wee (2000)} in turn uses the deepen
theory of analytic functions, including an untraditional result of
the theory of power series given by a theorem of Wiener.

Ramalhoto and Gomez-Corral \cite{Ramalhoto and Gomez-Corral
(1998)} discuss retrials in $M/M/r/d$ loss queues and present an
appropriate decomposition formulae for losses and delay in those
queueing systems. In the case of Markovian system with many
servers the results of Ramalhoto and Gomez-Corral \cite{Ramalhoto
and Gomez-Corral (1998)} are useful for analysis of the effect of
retrials in asymptotic analysis of losses.

The present paper provides the asymptotic analysis of the loss
probability of the $GI/M/m/n$ queue as $n\to\infty$, by reduction
of the sequence $\pi_{m,0,n}^{(n)}$ to above representation
(\ref{r1.5}), and then estimates the loss probability
$p_{m,n}=1/\pi_{m,0,n}^{(n)}$ by using the results of Tak\'acs
\cite{Takacs (1967)} and Postnikov \cite{Postnikov (1979)} given
in the Appendix. (The precise definition of the sequence
$\pi_{m,0,n}^{(n)}$ is given later.)

This is the same idea as in the earlier paper of Abramov
\cite{Abramov (2002)} related to the $GI/M/1/n$ queue, however it
is necessary to underline the following. Whereas in the case of
$GI/M/1/n$ the reduction to (\ref{r1.3}) is straightforward, and
the representation of the probabilities $r_n$ and their generating
function is very simple, the reduction to the recurrence relation
of (\ref{r1.5}) in the case of $GI/M/m/n$ queue, $m>1$, is not
obvious, and representation for probabilities $r_{k,m-k,n}$ and
$r_{0,m,j}$ and associated generating function is more difficult.
(The aforementioned probabilities $r_{k,m-k,n}$ and $r_{0,m,j}$
are defined later.) Furthermore, the sequence $\pi_{m,0,j}^{(n)}$,
$j\leq n$, is structured as schema with the series classes, and,
as $j=n$, the value $p_{m,n}=1/\pi_{m,0,n}^{(n)}$ coincides with
the desired loss probability. For this reason the class of
asymptotic results, that we could obtain here for $GI/M/m/n$
queue, is poorer than that for the $GI/M/1/n$ queue in Abramov
\cite{Abramov (2002)}.

Our approach has the following two essential advantages compared
to the pure analytical approaches of Choi et al \cite{Choi Kim and
Wee (2000)}, Kim and Choi \cite{Kim Choi (2003)}, Choi et al
\cite{Choi Kim Kim and Wee (2003)}, Simonot \cite{Simonot (1998)}:

\smallskip
$\bullet$ The problem reduces to the known classic results
(Theorem of Tak\'acs \cite{Takacs (1967)} and Tauberian Theorem of
Postnikov \cite{Postnikov (1979)}), permitting us to substantially
diminish the cumbersome algebraic calculations and clearly
understand the results.

\smallskip
$\bullet$ Along with standard asymptotic results having a
quantitative feature we also prove one interesting property
related to the case of $n$ increasing to infinity and the load
approaching 1 from the left. (For the more precise assumptions see
formulation of Theorem \ref{t2}.) Specifically it is proved that
the obtained asymptotic representation is the same for all $m\geq
1$, i.e. it coincides with asymptotic representation obtained
earlier for the $GI/M/1/n$ queueing system in Abramov
\cite{Abramov (2002)}. As $n$ increases to infinity, this
asymptotic property remains true for all fixed $\varrho\geq 1$.

\smallskip
The conditions of Theorem \ref{t2} and the first two cases of
Theorem \ref{t1} all fall into the domain of heavy traffic theory
(e.g. Borovkov \cite{Borovkov (1976)}, Whitt \cite{Whitt (2004)},
\cite{Whitt (2005)}). Although the aforementioned recent results
of Whitt \cite{Whitt (2004)}, \cite{Whitt (2005)} are related to
more general models, they however do not cover the results of this
paper.

Our approach is based on asymptotic analysis of relations
\eqref{r3.4} and \eqref{r3.5} which is based on asymptotic
representation for the root of equation $z=\alpha(\mu m-\mu mz)$
(see formulation of Theorem \ref{t1}) as $\varrho$ approaches 1,
presented in the book of Subhankulov \cite{Subhankulov (1976)}.
(Chapter 9 of this book is devoted to application of Tauberian
theorems to a specific moving server problem arising in Operations
Research.)

\smallskip

The rest of the paper is organized as follows. In Section 2 we
give some heuristic arguments preparing the reader to the results
of the paper. There are two theorems presenting the main results
in Section 3. In Section 4 we derive the recurrence relation of
convolution type for the loss probability. These recurrence
relations are then used to prove Theorem \ref{t1} of the paper.
Theorem \ref{t1} is proved in Section 5 and Section 6. In Section
7 we study the behavior of the loss probability as the load
approaches 1 from the left. In Section 8 numerical example
supporting the theory is provided. The Appendix contains auxiliary
results necessary for the purpose of the paper: the theorem of
Tak\'acs \cite{Takacs (1967)} and the Tauberian theorem of
Postnikov \cite{Postnikov (1979)}.

\section{\bf Some heuristic arguments}
\noindent The $GI/M/m/n$ queueing system, $m>1$, is more
complicated than its analog with one server. Letting $n$ be
infinity, discuss first the $GI/M/1$ and $GI/M/m$ queueing
systems.

It is well-known that the stationary queue-length distribution of
$GI/M/1$ queue (immediately before arrival of a customer with
large order number) is geometric. The same (customer-stationary)
queue-length distribution of $GI/M/m$ queue, provided that
immediately before arrival at least $m-1$ servers are occupied, is
geometric as well. Thus, a typical behavior of the queue-length
processes of $GI/M/1$ and $GI/M/m$ queues is similar, if we assume
additionally that a customer arriving into $GI/M/m$ queueing
system finds at least $m-1$ servers busy (see Kleinrock
\cite{Kleinrock (1975)}).

The similar situation holds in the case of the $GI/M/1/n$ and
$GI/M/m/n$ queues. Specifically, in the case of the $GI/M/1/n$
queue the stationary loss probability satisfies
(\ref{r1.1})-(\ref{r1.3}). In the case of the $GI/M/m/n$ queue the
{\it conditional} stationary loss probability {\it provided that
upon arrival at least $m-1$ servers are busy} satisfies
(\ref{r1.3}) as well. In the case of the $GI/M/m/n$ queue the only
difference is that, the value $\mu$ in (\ref{r1.2}) and
(\ref{r1.4}) should be replaced with $\mu m$, and $\sigma$ should
be the least in absolute value root of the equation $z=\alpha(\mu
m-\mu mz)$ rather than of the equation $z=\alpha(\mu-\mu z)$. In
the sequel, the least root of the functional equation
$z=\alpha(\mu m-\mu mz)$ is denoted $\sigma_m$.

Let us now discuss the stationary probabilities of the $GI/M/m/n$
system, $m>1$, without the condition above. It is clear that
eliminating the condition above proportionally changes the
stationary probabilities $\mathbf{P}\{$ arriving customer meets
$m+j-1$ customers in the system  \}, $j\geq 0$. That is, the loss
probability is changed proportionally as well. This enables us to
anticipate the behavior of the loss probability as $n\to\infty$ in
some cases.

Specifically, in the case $\varrho<1$ and $n$ large, the loss
probability is equal to conditional loss probability, provided
that upon arrival at least $m-1$ servers are busy, multiplied by
some constant. That is, as $n$ large, the both abovementioned loss
probabilities, conditional and unconditional, are of the same
order. Following Abramov \cite{Abramov (2002)} (see also Choi et
al \cite{Choi Kim Kim and Wee (2003)}) this order is
$O(\sigma_m^n)$. The precise result is given by Theorem \ref{t1}
below.

If $n$ large and $\varrho\geq 1$, then the probability that
arriving customer meets less than $m-1$ customers in the system is
small, and therefore, the loss probability should be approximately
the same as the conditional stationary probability that upon
arrival of a customer at least $m-1$ servers are busy. That is,
following Abramov \cite{Abramov (2002)} (see also Choi et al
\cite{Choi Kim Kim and Wee (2003)}) one can expect, that in the
case of $\varrho\geq1$, the limiting stationary loss probability
of the $GI/M/m/n$ queue, as $n\to\infty$, should be equal to
$(\varrho-1)/\varrho$ for all $m$. Can the last property of
asymptotic independence of $m$ be extended as $n$ increases to
infinity and $\varrho$ approaches 1? The paper provides the
condition for this asymptotic independence in this case.

\section{\bf Formulation of the main results}

\begin{thm}\label{t1}
 If $\varrho>1$ then for any $m\geq1$
\begin{equation}\label{r3.1}
 \lim_{n\to\infty}p_{m,n}=\frac{\varrho-1}{\varrho}.
\end{equation}
If $\varrho=1$ and $\varrho_2=\int_0^\infty (\mu
x)^2\mbox{d}A(x)<\infty$  then for any $m\geq1$
\begin{equation}\label{r3.2}
\lim_{n\to\infty}np_{m,n}=\frac{\varrho_2}{2}.
\end{equation}
If $\varrho=1$ and $\varrho_3=\int_0^\infty (\mu
x)^3\mbox{d}A(x)<\infty$ then for large $n$ and any $m\geq1$ we
have
\begin{equation}\label{r3.3}
p_{m,n}=\frac{\varrho_2}{2n}+O\Big(\frac{\log n}{n^2}\Big).
\end{equation}
If $\varrho<1$ then for $p_{m,n}$ we have the limiting relation:
\begin{equation}\label{r3.4}
\lim_{n\to\infty}\frac{p_{m,n}}{\sigma_m^n}=K_m[1+\mu m\alpha'(\mu
m-\mu m\sigma_m)],
\end{equation}
where $\alpha^\prime(\cdot)$ denotes the derivative of
$\alpha(\cdot)$,
\begin{equation}\label{r3.5}
K_m=\Big[1+(1-\sigma_m)\sum_{j=1}^m\frac{
\binom{m}{j}C_j}{(1-\varphi_{j})}~\frac{m(1-\varphi_{j})-j}{m(1-\sigma_m)-j}\Big]^{-1},
\end{equation}
\[
\varphi_{j}=\int_0^\infty\mbox{e}^{-\mu jx}\mbox{d}A(x),
\]
\[
C_j=\prod_{i=1}^j\frac{1-\varphi_{j}}{\varphi_{j}},
\]
and $\sigma_m$ is the least in absolute value root of functional
equation:
\[
z=\alpha(\mu m-\mu mz).
\]
\end{thm}

Theorem \ref{t1} shows that if $\varrho>1$ then the limiting
stationary loss probability is independent of parameter $m$. If
$\varrho=1$ and $\varrho_2<\infty$ then
$\lim_{n\to\infty}np_{m,n}$ is independent of parameter $m$ as
well. The proof of (\ref{r3.1}) seems to be given by simple
straightforward arguments (extended version of the heuristic
arguments of Section 2). Nevertheless, all results are proved by
reduction to the abovementioned theorems of Tak\'acs \cite{Takacs
(1967)} and Postnikov \cite{Postnikov (1979)} given in the
Appendix. The most significant result of Theorem \ref{t1} is
(\ref{r3.4}). This result is then used to prove the statements of
Theorem \ref{t2} on the behavior of the loss probability as the
load approaches 1 from the left.

This behavior of the loss probability is given by the following
theorem.

\begin{thm}\label{t2}
 Let $\varrho=1-\varepsilon$, where
$\varepsilon>0$, and $\varepsilon n\to C$ as $n\to\infty$ and
$\varepsilon\to 0$. Assume that $\varrho_3=\varrho_3(n)$ is a
bounded sequence in $n$, and there exists $\widetilde
\varrho_2=\lim_{n\to\infty}\varrho_2(n)$. In the case where $C>0$
for any $m\geq1$ we have
\begin{equation}\label{r3.9}
p_{m,n}=\frac{\varepsilon\mbox{e}^{-2C/\widetilde\varrho_2}}
{1-\mbox{e}^{-2C/\widetilde\varrho_2}}[1+o(1)].
\end{equation}
In the case where $C=0$ for any $m\geq1$ we have
\begin{equation}\label{r3.10}
p_{m,n}=\frac{\widetilde\varrho_2}{2n}+o\Big(\frac{1}{n}\Big).
\end{equation}

\end{thm}

Theorem \ref{t2} shows that as $\varrho$ approaches  1 from the
left, the loss probability $p_{m,n}$  becomes independent of
parameter $m$ when $n$ large, and the asymptotic behavior of the
loss probability is exactly the same as for the $GI/M/1/n$ queue.

\section{\bf Derivation of the recurrence equations for the loss
probability} \noindent For the sake of convenience, in this
section we keep in mind that the first $m-1$ states of the
$GI/M/m/n$ queue-length process form one special class. If an
arriving customer occupies one of servers, then the system is
assumed to be in this class, and the states of this class are
numbered 1, 2,..., $m$. Otherwise, the system is in the other
class with states $m+1$, $m+2$,..., $m+n$, where the last state,
$m+n$, is associated with a loss of an arriving customer.

For example, if $n=0$, then the second class of the $GI/M/m/0$
queueing system consists of one state only.

Let us now build the recurrence relation similar to those of the
$GI/M/1/n$ queue.

We start from the $GI/M/1/0$ queue. For this queue we have
\begin{equation}\label{r4.1}
 \pi_{1,0}=\frac{1}{r_{0,1}},
\end{equation}
where
\[
r_{0,1}=\varphi_1=\int_0^\infty\mbox{e}^{-\mu x}\mbox{d}A(x).
\]
Equation (\ref{r4.1}) formally follows from the recurrence
relation $ \pi_{0,1}=r_{0,1}\pi_{1,0} $, where $\pi_{0,1}=1$. The
loss probability for the $GI/M/1/0$ queue is equal to
\[
p_{1,0}=\frac{1}{\pi_{1,0}}=\varphi_1.
\]

Before considering the case of the $GI/M/m/n$ queue, notice that
the value $\pi_{m-k,k}$ has the meaning of the expected number of
arrivals to the stationary system up to at the first time an
arriving customer  finds $m-k$ servers busy and $k$ remaining
servers free.

In the case of the $GI/M/2/0$ queue, by the total expectation
formula we have
\[
\pi_{1,1}=r_{0,2}\pi_{2,0}+r_{1,1}\pi_{1,1},
\]
where
\[
r_{1,1}=2\int_0^\infty[1-\mbox{e}^{-\mu x}]\mbox{e}^{-\mu
x}\mbox{d}A(x),
\]
\[
r_{0,2}=\varphi_2=\int_0^\infty\mbox{e}^{-2\mu x}\mbox{d}A(x).
\]

Then, by the total expectation formula, the recurrence relation
for the $GI/M/m/0$ queue looks as follows:
\begin{equation}\label{r4.7}
 \pi_{m-1,1}=\sum_{k=0}^{m-1}r_{k,m-k}\pi_{m-k,k},
\end{equation}
where
\[
%P_{i,j}
r_{k,m-k}= \binom{m}{k}\int_0^\infty[1-\mbox{e}^{-\mu
x}]^{k}\mbox{e}^{-(m-k)\mu x}\mbox{d}A(x).
\]

It is well-known that
\begin{equation}\label{r4.9}
\pi_{m,0}= \sum_{i=0}^m\binom{m}{i} \prod_{j=1}^i
\frac{1-r_{0,j}}{r_{0,j}}=\sum_{i=0}^m\binom{m}{i}
 C_i,
\end{equation}
and the loss probability is
\begin{equation}\label{r4.10}
p_{m,0}=\frac{1}{\pi_{m,0}}=\Big[\sum_{i=0}^m\binom{m}{i}
 C_i\Big]^{-1}
\end{equation}
(see  Cohen \cite{Cohen (1957)}, Palm \cite{Palm (1943)},
Pollaczek \cite{Pollaczek (1953)}, Tak\'acs \cite{Takacs (1957)}
as well as Bharucha-Reid \cite{Bharucha-Reid (1960)}).

A relatively simple proof of (\ref{r4.9}) and (\ref{r4.10}) can be
found in Tak\'acs \cite{Takacs (1957)}. It is based on another
representation than (\ref{r4.7}). For our further purposes,
representation (\ref{r4.7}) is preferable. Representation
(\ref{r4.7}) is a recurrence relation of the convolution type
(\ref{r1.5}), and in the following it helps us to reduce the
problem to the abovementioned combinatorial results of Tak\'acs
\cite{Takacs (1967)}. Once this is done, we apply then the
Tauberian theorem of Postnikov \cite{Postnikov (1979)}.

\medskip

Let us now consider the $GI/M/m/n$ queueing system. In the case of
this system with $n\geq 1$ we add an additional subscript to the
notation. Specifically, $r_{k,m-k,0}=r_{k,m-k}$, and for $j\leq n$
\[
r_{0,m,j}=\int_0^\infty \mbox{e}^{-m\mu x}\frac{(m\mu
x)^j}{j!}\mbox{d}A(x),
\]
and
\[
r_{k,m-k,n} =\binom{m}{k}\int_0^\infty\mbox{e}^{-(m-k)\mu x}
\]
\[\times
\Big\{\int_0^x\frac{(m\mu u)^{n-1}}{(n-1)!}(\mbox{e}^{-\mu
u}-\mbox{e}^{-\mu x})^{k}m\mu\mbox{d}u\Big\}\mbox{d}A(x).
\]
In addition, the value $\pi_{m,0,n}^{(n)}$ denotes the expected
number of arrivals into the stationary $GI/M/m/n$ queue up to the
first loss. (Replacing the indexes $n$ with $j$ has the same
meaning for the $GI/M/m/j$ queue.) Also, we will use the notation
$\pi_{m-k,k,0}^{(0)}$ instead of the earlier notation
$\pi_{m-k,k}$ for the $GI/M/m/0$ queue.

Then the recurrence relation associated with the $n$th series
looks as follows:
\begin{equation}\label{r4.13}
\pi_{m,0,j}^{(n)}=\sum_{l=0}^{j}r_{0,m,l}\pi_{m,0,j-l+1}^{(n)}+\sum_{k=1}^{m-1}r_{k,m-k,j}\pi_{m-k,k,0}^{(n)},
\end{equation}
\[
j=0,1,...,n-1,
\]
where
\begin{equation}\label{r4.14}
\pi_{m-i,i,0}^{(n)}=\sum_{k=0}^{m-i}r_{k,m-k-i+1,0}\pi_{m-k-i+1,k+i-1}^{(n)}
\ \ (\pi_{0,m,0}^{(n)}=1),
\end{equation}
\[
i=1,2,...,m-1,
\]
and the second sum of (\ref{r4.13}) is equal to zero if $m=1$.
Moreover, if $m=1$, then we do not longer need the upper index
$(n)$, showing the series number, and equation (\ref{r4.14}). It
is not difficult to see, that for the given series $n$, the
recurrence relations (\ref{r4.13}) and (\ref{r4.14}) form a
recurrence relation of the convolution type given by (\ref{r1.5}).
In the next section we prove relations (\ref{r3.1})-(\ref{r3.3})
of Theorem \ref{t1}.

\section{\bf The proof of (\ref{r3.1})-(\ref{r3.3})}
\noindent
First of all note, that
\[
\lim_{n\to\infty}\Big[\sum_{l=0}^{n}r_{0,m,l}+\sum_{k=1}^{m-1}r_{k,m-k,n}\Big]
= \sum_{l=0}^{\infty}r_{0,m,l}
\]
\[=\sum_{l=0}^\infty\int_0^\infty\mbox{e}^{-m\mu
x}\frac{(m\mu x)^l}{l!}\mbox{d}A(x)
\]
\[=
\int_0^\infty\sum_{l=0}^\infty\mbox{e}^{-m\mu x}\frac{(m\mu
x)^l}{l!}\mbox{d}A(x)=1.
\]
Therefore, one can apply the theorem of Tak\'acs \cite{Takacs
(1967)} (see Appendix). Let $\gamma_1$ denote
\begin{equation}\label{r5.2}
\gamma_1=\sum_{l=1}^\infty lr_{0,m,l}.
\end{equation}
Then also
\begin{equation}\label{r5.3}
\lim_{n\to\infty}\Big[\sum_{l=1}^{n}lr_{0,m,l}+\sum_{k=1}^{m-1}(n+k)r_{k,m-k,n}\Big]=\gamma_1.
\end{equation}
This is because
\[
(n+k)\binom{m}{k}\int_0^\infty\mbox{e}^{-(m-k)\mu
x}\Big\{\int_0^x\frac{(m\mu u)^{n-1}}{(n-1)!}(\mbox{e}^{-\mu
u}-\mbox{e}^{-\mu x})^{k}m\mu\mbox{d}u\Big\}\mbox{d}A(x)
\]
\[\leq
(n+k)\binom{m}{k}\int_0^\infty\mbox{e}^{-\mu
(m-k)x}\Big\{\int_0^x\frac{(m\mu
u)^{n-1}}{(n-1)!}m\mu\mbox{d}u\Big\}\mbox{d}A(x)
\]
\[=
(n+k)\binom{m}{k}\int_0^\infty\mbox{e}^{-\mu (m-k)x}\frac{(m\mu
x)^{n}}{n!}\mbox{d}A(x)\to 0,
\]
as $n\to\infty$.

According to (\ref{r5.2}) and (\ref{r5.3}) we have
$\gamma_1=m\mu/\lambda$, and therefore, $\gamma_1=1/\varrho$. Then
according to theorem of Tak\'acs \cite{Takacs (1967)}, given in
the Appendix, in the case of $\varrho>1$ we obtain
\[
\lim_{n\to\infty}\pi_{m,0,n}^{(n)}=\frac{1}{1-\gamma_1}=\frac{\varrho}{\varrho-1}.
\]
Then, in this case of the limiting loss probability as
$n\to\infty$ we obtain
\[
\lim_{n\to\infty}p_{m,n}=\lim_{n\to\infty}\frac{1}{\pi_{m,n,0}^{(n)}}=\frac{\varrho-1}{\varrho}.
\]
Similarly to (\ref{r5.2}), Let $\gamma_2$ denote
\[
\gamma_2=\sum_{l=2}^\infty l(l-1)r_{0,m,l}.
\]
Then also
\[
\lim_{n\to\infty}\Big[\sum_{l=2}^{n}l(l-1)r_{0,m,l}+\sum_{k=1}^{m-1}(n+k)(n+k-1)r_{k,m-k,n}\Big]=\gamma_2,
\]
and $\varrho_2<\infty$ as $\gamma_2<\infty$. Indeed, as
$n\to\infty$,
\[
(n+k)(n+k-1)
\]
\[\times~
\binom{m}{k}\int_0^\infty\mbox{e}^{-(m-k)\mu
x}\Big\{\int_0^x\frac{(m\mu u)^{n-1}}{(n-1)!}(\mbox{e}^{-\mu
u}-\mbox{e}^{-\mu x})^{k}m\mu\mbox{d}u\Big\}\mbox{d}A(x)
\]
\[\leq
(n+k)(n+k-1)\binom{m}{k}\int_0^\infty\mbox{e}^{-\mu
(m-k)x}\frac{(m\mu x)^{n}}{n!}\mbox{d}A(x)\to 0.
\]
Therefore, in the case where $\varrho=1$ and $
\varrho_2=\int_0^\infty (\mu mx)^2\mbox{d}A(x)<\infty $ we obtain
\begin{equation}\label{r5.10}
\lim_{n\to\infty}np_{m,n}=\frac{\varrho_2}{2}.
\end{equation}
Limiting relation (\ref{r5.10}) can be improved with the aid of
the Tauberian theorem of Postnikov \cite{Postnikov (1979)} (see
Appendix). In the case where $\varrho=1$ and $
\varrho_3=\int_0^\infty (\mu mx)^3\mbox{d}A(x)<\infty$, for large
$n$ we obtain
\[
p_{m,n}=\frac{\varrho_2}{2n}+O\Big(\frac{\log n}{n^2}\Big).
\]
Indeed, let $\gamma_3$ denote
\[
\gamma_3=\sum_{l=3}^\infty l(l-1)(l-2)r_{0,m,l}.
\]
Then also
\[
\lim_{n\to\infty}\Big[\sum_{l=3}^{n}l(l-1)(l-2)r_{0,m,l}
\]
\[
+\sum_{k=1}^{m-1}(n+k)(n+k-1)(n+k-2)r_{k,m-k,n}\Big]=\gamma_3,
\]
and $\varrho_3<\infty$ as $n\to\infty$. Indeed, as $n\to\infty$,
\[
(n+k)(n+k-1)(n+k-2)
\]
\[\times~
\binom{m}{k}\int_0^\infty\mbox{e}^{-(m-k)\mu
x}\Big\{\int_0^x\frac{(m\mu u)^{n-1}}{(n-1)!}(\mbox{e}^{-\mu
u}-\mbox{e}^{-\mu x})^{k}m\mu\mbox{d}u\Big\}\mbox{d}A(x)
\]
\[\leq
(n+k)(n+k-1)(n+k-2)\binom{m}{k}
\]
\[\times \int_0^\infty\mbox{e}^{-\mu (m-k)x}\frac{(m\mu
x)^{n}}{n!}\mbox{d}A(x)\to 0.
\]
Thus, (\ref{r3.2}) and (\ref{r3.3}) follow.

\section{\bf The proof of (3.4)}
\noindent Whereas (\ref{r3.1})-(\ref{r3.3}) are proved by
immediate reduction to the known results associated with
(\ref{r1.5}), the proof of (\ref{r3.4}) requires special analysis.
In order to simplify the analysis let us concentrate our attention
to the constant $K_m$ in relation (\ref{r3.4}). Multiplying this
constant by $(1-\sigma_m)$ we obtain
\[ \widetilde K_m=(1-\sigma_m)K_m
\]
\begin{equation}\label{r6.1}
=
\Big[\frac{1}{1-\sigma_m}+\sum_{j=1}^m\frac{\binom{m}{j}
C_j}{(1-\varphi_{j})}~\frac{m(1-\varphi_{j})-j}{m(1-\sigma_m)-j}\Big]^{-1}.
\end{equation}
 The constant $\widetilde K_m$, given by (\ref{r6.1}),
is well-known from the theory of $GI/M/m$ queueing system.
Specifically, let $\widetilde p_j$, $j=0,1,...,$ be the stationary
probabilities of the number of customers in this system
immediately before arrival of a customer. It is known (e.g.
Bharucha-Reid \cite{Bharucha-Reid (1960)}, Borovkov \cite{Borovkov
(1976)}) that for all $j\geq m$
\begin{equation}\label{r6.2}
\widetilde p_j=\widetilde K_m\sigma_m^{j-m}.
\end{equation}

Now, in order to prove (\ref{r3.4}) let us write a new recurrence
relation, alternative to (\ref{r4.13}). For this purpose, join the
first $m$ states of the $GI/M/m/n$ process to a single state and
label it 0. Other states will be numbered 1, 2,..., $n$. In the
new terms we have the following recurrence relations
\begin{equation}\label{r6.3}
\Pi_j^{(n)}=\sum_{i=0}^{j} r_{0,m,i}\Pi_{j-i+1}^{(n)},
\end{equation}
with some initial value $\Pi_0^{(n)}$ for the given series $n$.
For example, for the series $n=0$ we have
\[
\Pi_0^{(0)}=\sum_{i=0}^m\binom{m}{i}
 C_i
\]
(see relation (\ref{r4.9})). A formal application of the theorem
of Tak\'acs \cite{Takacs (1967)} (see Appendix), applied to
recurrence relation (\ref{r6.3}), for large $n$ yields:
\begin{equation}\label{r6.5}
\lim_{n\to\infty}\frac{\Pi_n^{(n)}\sigma_m^n}{\Pi_0^{(n)}}=\frac{1}{1+\mu
m\alpha'(\mu m-\mu m\sigma_m)}.
\end{equation}

Let us now find $\lim_{n\to\infty}\Pi_0^{(n)}$. Notice, that for
$j\geq m$ the probability $\widetilde p_j$ can be rewritten as
follows. From (\ref{r6.2}) we have:
\begin{equation}\label{r6.6}
\widetilde p_j = \widetilde
K_m\sigma_m^{j-m}=K_m\sigma_m^{j-m}(1-\sigma_m) =
\frac{K_mP_j}{\sigma_m^m},
\end{equation}
where $P_j$ is the conditional probability for the $GI/M/m$ queue,
that an arriving customer finds $j$ customers in the queue
provided that upon arrival at least $m-1$ servers are occupied.
The conditional probability $P_{j-m}$ coincides with the
stationary queue-length distribution immediately before arrival of
a customer in the $GI/M/1$ queue given under the expected service
time $(\mu m)^{-1}$. $K_m$ is the stationary probability for the
$GI/M/m$ queue, that upon arrival at least $m-1$ servers are
occupied.

From the theory of Markov chains associated with $GI/M/1/n$ queue
it is known (e.g. Choi and Kim \cite{Choi and Kim (2000)}) that
the $j$-state probability immediately before arrival of a customer
is $(\pi_{n-j}-\pi_{n-j-1})/\pi_n$, where $\pi_n$ is given by
(\ref{r1.3}), and in turn the loss probability is determined by
(\ref{r1.1}).

Then for the same $j$-state probability of $GI/M/1$ queue with the
expected service time $(\mu m)^{-1}$ we have
\begin{equation}\label{r6.7}
P_j=
\lim_{n\to\infty}\frac{\Pi_{n-j}^{(n)}-\Pi_{n-j-1}^{(n)}}{\Pi_{n}^{(n)}}.
\end{equation}

In turn, by (\ref{r6.6}) and (\ref{r6.7}), the $j+m$-state
probability of the $GI/M/m$ queue is determined as follows:
\[
\widetilde p_{j+m}=K_m\lim_{n\to\infty}\frac{\Pi_{n-j}^{(n)}
-\Pi_{n-j-1}^{(n)}}{\Pi_{n}^{(n)}},
\]
and
\begin{equation}\label{r6.9}
\lim_{n\to\infty}\Pi_0^{(n)}=\frac{1}{K_m}.
\end{equation}
In view of (\ref{r6.9}) and (\ref{r6.5}) and according to Tak\'acs
theorem \cite{Takacs (1967)} we obtain:
\begin{equation}\label{r6.10}
\lim_{n\to\infty}\left[\Pi_n^{(n)}-\frac{1}{K_m\sigma_m^n[1+\mu
m\alpha'(\mu m-\mu m\sigma_m)]}\right]=\frac{\varrho
K_m}{1-\varrho}.
\end{equation}
Now, taking into consideration that the loss probability
\[
p_{m,n}=\frac{1}{\Pi_n^{(n)}},
\]
we obtain statement (\ref{r3.4}) of the theorem.

\section{\bf The proof of Theorem \ref{t2}}
\noindent It was shown in Subhankulov \cite{Subhankulov (1976)},
p. 326, that if $\varrho^{-1}=1+\varepsilon$, $\varepsilon>0$ and
$\varepsilon\to 0$, $\varrho_3(n)$ is a bounded sequence, and
there exists $\widetilde \varrho_2=\lim_{n\to\infty}\varrho_2(n)$,
then
\begin{equation}\label{r7.1}
\sigma_m=1-\frac{2\varepsilon}{\widetilde\varrho_2}+O(\varepsilon^2),
\end{equation}
where $\sigma_m=\sigma_m(n)$ is the minimum in absolute value root
of the functional equation $z=\alpha(\mu m-\mu mz)$, $|z|\leq 1$,
and where the parameter $\mu$ and the function $\alpha(z)$, both
or one of them, are assumed to depend on $n$. (Asymptotic
representation \eqref{r7.1} can be immediately obtained by
expanding the equation $z-\alpha(\mu m-\mu mz)=0$ for small $z$.)

Then, after some algebra we have
\begin{equation}\label{r7.2}
[1+\mu m\alpha'(\mu m-\mu m\sigma_m)]=\varepsilon+o(\varepsilon),
\end{equation}
and
\begin{equation}\label{r7.3}
\sigma_m^n=\mbox{e}^{-2C/\varrho_2}[1+o(1)].
\end{equation}
In view of (\ref{r7.1}), the term
\[
(1-\sigma_m)\sum_{j=1}^m\frac{\binom{m}{j}C_j}{(1-\varphi_{j})}~\frac{m(1-\varphi_{j})-j}{m(1-\sigma_m)-j}
\]
has the order $O(\varepsilon)$. Therefore, for term (\ref{r3.5})
we have
\begin{equation}\label{r7.4}
K_m=1+O(\varepsilon),
\end{equation}
and in the case where $C>0$, in view of (\ref{r7.2})-(\ref{r7.4})
and \eqref{r6.10}, we obtain:
\begin{equation}\label{r7.5}
p_{m,n}=\frac{\varepsilon\mbox{e}^{-2C/\widetilde\varrho_2}}
{1-\mbox{e}^{-2C/\widetilde\varrho_2}}[1+o(1)].
\end{equation}
(\ref{r3.9}) is proved.

The proof of (\ref{r3.10}) follows by expanding the main term of
the asymptotic expression of (\ref{r7.5}) for small $C$.

Theorem \ref{t2} is completely proved.

\section{\bf Numerical example}
\noindent In this section a numerical example supporting the
theory is provided. Specifically, we simulate $D/M/1/n$ and
$D/M/2/n$ queues and check statements (\ref{r3.10}) and
(\ref{r3.9}) of Theorem \ref{t2} numerically.  The results of
simulation are reflected in the table below. The value $\varrho$
is taken 0.999, so that $\epsilon=0.001$ The value $n$ varies from
10 to 50, and parameter $C=\epsilon n$ varies from 0.01 to 0.05
The theoretical values of the loss probability for these $n$ are
calculated by (\ref{r3.10}). There are also the loss probabilities
for $n=100$. The theoretical value for the loss probability
related to this case is calculated by (\ref{r3.9}).

\begin{table}
    \begin{center}
        \begin{tabular}{c||c|c|c}\hline
& Loss probability & Loss probability   & Loss probability\\
$n$& theoretical & simulated for  & simulated for \\
& & $D/M/1/n$ queue &$D/M/2/n$ queue\\
\hline
10        & 0.0501   & 0.0426  & 0.0390  \\
15        & 0.0334   & 0.0292  & 0.0275  \\
20        & 0.0251   & 0.0221  & 0.0211  \\
25        & 0.0200   & 0.0180  & 0.0173  \\
30        & 0.0167   & 0.0151  & 0.0146  \\
35        & 0.0143   & 0.0128  & 0.0124  \\
40        & 0.0125   & 0.0111  & 0.0108  \\
45        & 0.0111   & 0.0098  & 0.0096  \\
50        & 0.0100   & 0.0087  & 0.0085  \\
100       & 0.0045   & 0.0040  & 0.0039  \\
\hline
        \end{tabular}

        \caption{The comparison table of the loss probabilities for $D/M/1/n$ and
        $D/M/2/n$ queues}
    \end{center}
\end{table}

The table is structured as follows: Column 1 contains the values
of parameter $n$, Column 2 contains the theoretical values for the
loss probability given by (\ref{r3.10}) and (\ref{r3.9}), Column 3
and 4 contain the loss probabilities obtained by simulation for
the $D/M/1/n$ and $D/M/2/n$ queueing systems respectively.

As we can see from this table the difference between the loss
probabilities of the single-server and two-server queueing systems
obtained by simulation is not large, and difference between these
loss probabilities decreases as $n$ increases. As $n$ increases
the both simulated loss probabilities approach the theoretical
loss probability.

\section*{\bf Acknowledgements}
The author thanks the referees for useful comments. Especial thank
is to the referee calling attention of the author to the results
of Whitt \cite{Whitt (2004)}, \cite{Whitt (2005)} having immediate
relation to the main result of this paper. The initiation of the
paper was due to a question (conjecture) of Professor Henk C.
Tijms related to asymptotic behavior of the loss probability in
the overloaded $GI/M/m/n$ queue. The question had relation to the
talk of the author at the First Madrid Conference on Queueing
Theory.

\section*{\bf APPENDIX}

In the appendix we recall the main results on asymptotic behavior
of the sequence $Q_n$, as $n\to\infty$ (see relation
(\ref{r1.5})).

\medskip

Denote $f(z)=\sum_{i=0}^\infty f_iz^i, \ \ |z|\leq 1,$
$\gamma_i=\sum_{j=i}^\infty\Big(\prod_{k=j-i+1}^{j}k\Big)f_j.$

\noindent {\bf Theorem A1.} (Tak\'acs \cite{Takacs (1967)}, p. 22,
23.) \textit{If $\gamma_1<1$ then
$$
\lim_{n\to\infty}Q_n=\frac{Q_0}{1-\gamma_1}.
$$
If $\gamma_1=1$ and $\gamma_2<\infty$, then
$$
\lim_{n\to\infty}\frac{Q_n}{n}=\frac{2Q_0}{\gamma_2}.
$$
If $\gamma_1>1$ then
$$
\lim_{n\to\infty}\left[Q_n-\frac{Q_0}{\delta^n(1-f'(\delta))}\right]=\frac{Q_0}{1-\gamma_1},
$$
where $\delta$ is the least in absolute value root of equation
$z=f(z)$.}

\medskip

\noindent {\bf Theorem A2.} (Postnikov \cite{Postnikov (1979)},
Section 25.) \textit{If $\gamma_1=1$ and $\gamma_3<\infty$, then
as $n\to\infty$}
$$
Q_n=\frac{2Q_0}{\gamma_2}n+O(\log n).
$$

%\section*{Acknowledgement}

%An initiation of this paper was due to the conjecture of Professor
%Henk Tijms (Vrije University) that in the case of $\varrho>1$, the
%limiting stationary loss probability of the $GI/M/m/n$ queue, as
%$n\to\infty$, should be equal to the same value
%$(\varrho-1)/\varrho$ for all $m$.

\bibliographystyle{amsplain}

\newpage

\section*{\bf Bibliography}

\textbf{Vyacheslav M. Abramov} graduated from Tadzhik State
University (Dushanbe, Tadzhikistan) in 1977. During the period
1977-1992 he worked at the Research Institute of Economics under
the Tadzhikistan State Planning Committee (GosPlan). In 1992 he
repatriated to Israel and during 1994-2001 worked in software
companies of Israel as a software engineer and algorithms
developer. In 2002-2005 he was an assistant and  lecturer in Judea
and Samaria College, Tel Aviv University and Holon Institute of
Technology. In 2004 he received a PhD degree from Tel Aviv
University, and since 2005 has been working at School of
Mathematical Sciences of Monash University (Australia). The
scientific interests of him are mainly focused on the theory and
application of queueing systems. He is an author of a monograph
and various papers published in Journal of Applied Probability,
Annals of Operations Research, Queueing Systems, SIAM Journal on
Applied Mathematics and other journals.

\end{document}